\documentclass[11pt,a4paper]{article}
\usepackage{amsfonts,amsmath,amsthm,epsf,graphics,verbatim,amssymb,amscd,eucal,bbm}

\newtheorem{theorem}{Theorem}[section]

\newtheorem{proposition}[theorem]{Proposition}

\theoremstyle{definition}

\theoremstyle{remark}

\def\Blackboardfont{\mathbb}

\newcommand{\moins}{ {\setminus} }

\newcommand{\pres}[2]{\langle \: #1 \mid #2 \: \rangle}

\def\next{\text{Next}}

\def\Z{{\Blackboardfont Z}}

\def\N{{\Blackboardfont N}}
\def\R{{\Blackboardfont R}}

\def\cB{{\mathcal B}}
\def\cG{{\mathcal G}}

\def\cS{{\mathcal S}}

\def\cX{{\mathcal X}}

\def\cM{{\mathcal M}}

\def\cM{{\mathcal M}}

\def\first{\text{First}}

\def\eref#1{(\ref{#1})}

\def\norm#1{\|#1\|}

\begin{document}
\sloppy

\title{\bf Appendix to the paper ``Randomly Growing Braid on Three Strands and the Manta Ray''}

\author{Jean {\sc Mairesse}
\thanks{LIAFA, CNRS-Universit\'e Paris 7, case
    7014, 2, place Jussieu, 75251 Paris Cedex 05, France. E-mail: {\tt Jean.Mairesse@liafa.jussieu.fr}} 
    \ and Fr\'ed\'eric {\sc Math\'eus} 
\thanks{LMAM, Universit\'e de Bretagne-Sud, Campus de Tohannic,  BP
  573, 56017 Vannes, France. E-mail: {\tt Frederic.Matheus@univ-ubs.fr}}}

\maketitle


\begin{abstract}
This paper is an appendix to the paper ``Randomly Growing Braid on Three Strands and the Manta
  Ray'' by J. Mairesse and F. Math\'eus (to appear in the Annals of Applied Probability).
It contains the details of some computations, and the proofs of some
results concerning the examples treated there, as well as some
  extensions. 
\end{abstract}

\smallskip


\smallskip

\textsl{AMS classification (2000):} Primary 20F36, 20F69, 60B15;
Secondary 60J22, 82B41, 37M25.


\section{Introduction}
\label{se-intro}

Consider the braid group $B_3=\pres{a,b}{aba=bab}$ and the nearest
neighbor random walk defined by a probability
$\nu$ with support $S\: =\: \{a,a^{-1},b,b^{-1}\}$.
Let $(X_n)_n$ be a realization of the random walk. 
In order
to understand the asymptotic behavior of $X_n$, the first step is to
study its complexity, i.e. the length $|X_n|$ of $X_n$ with respect to
the generating system $S$. To this aim, a main quantity of interest is
the growth rate of the length: $\gamma=\lim_n |X_n|/n$. 

\par

In the paper ``Randomly Growing Braid on Three Strands and the Manta Ray'', we compute explicitly $\gamma$ for any
distribution $\nu$ on
$\{a,a^{-1},b,b^{-1}\}$. More precisely, given $\nu$, we define
eight polynomial equations of degree 2 over eight
indeterminates, see \cite[(30)]{MaManta}, which are shown to admit a
unique solution $r$. The
rate $\gamma$ is then obtained as an explicit linear functional of
$r$, see \cite[(29)]{MaManta}.

\par

Then we mention that the polynomial equations can be
completely solved to provide a closed form formula for $\gamma$ under
the following two natural symmetries:
(i) $\nu(a)=\nu(a^{-1})$ and $\nu(b)=\nu(b^{-1})$; (ii)
$\nu(a)=\nu(b)$ and $\nu(a^{-1})=\nu(b^{-1})$.
This is stated in \cite[Prop. 4.4]{MaManta} and \cite[Prop. 4.5]{MaManta}.
In the present Appendix, we give the proof of these two propositions, see
Sections \ref{se-prop44} and \ref{se-prop45}.

\par

The family of Artin groups of dihedral type $I_2(n)$ is a natural
generalization of the braid group $B_3$ and our techniques adapt to compute
the drift. The case of the simple random walk is stated in \cite[Prop. 5.4]{MaManta} and the proof is given in Section \ref{se-prop54}.

\par

The last Section contains additional informations concerning the
random walk on the quotient space $B_{3}/Z$ of $B_{3}$ by its center
$Z$. We compute the entropy, the minimal
positive harmonic functions, and the Green function. A central limit
theorem is also discussed.

This appendix is self-contained in the following sense: the examples
treated and the results proved are specified. On the other
hand, we make a
constant use of notations, notions or results from \cite{MaManta} which are neither
redefined nor reproved. Hence, it is certainly not a good idea to read this
appendix independently of the parent paper.

\section{Proof of Proposition 4.4} \label{se-prop44}

We prove the following statement which is Prop. 4.4 in \cite{MaManta}.
The drifts  $\gamma_{\Sigma}$, $\gamma_{\Delta}$ and $\gamma_{S_+}$ are defined in \cite[Eq. (14)]{MaManta}.

\medskip

{\bf Proposition 4.4.}
{\it Assume that $\nu(a)=\nu(a^{-1})=p$,
  $\nu(b)=\nu(b^{-1})=1/2-p$, with $p\in (0,1/4]$. Let $u$ be the
  smallest root in $(0,1)$ of the polynomial
\[
P=2(4p-1)X^3 +
  (24p^2-18p+1)X^2+ p(-12p+7)X + p(2p-1)\:.
\]
The various drifts are given by:
\begin{eqnarray*}
\gamma(p)=\gamma_{\Sigma}(p)=-2\gamma_{\Delta}(p)=\frac{2}{3}\gamma_{S_+}(p)=p+(1-4p)u\:.
\end{eqnarray*}
For $p\in [1/4,1/2)$, we have $\gamma(p)=\gamma(1/2-p)$ and similarly
  for $\gamma_{\Sigma}, \gamma_{S_+},$ and $\gamma_{\Delta}$.}

\begin{proof}
Assume that $\nu(a)=\nu(a^{-1})=p$ and $\nu(b)=\nu(b^{-1})=q=1/2-p$.
Let $r$ be the unique solution in
$\{x \in (\R_+^*)^{\Sigma} \mid  \sum_{u\in \Sigma} x(u) =1\}$ to
the Traffic Equations \cite[(30)]{MaManta}.
By symmetry, we should have:
\[
r(a)= r(ba\Delta), \quad r(b)=r(ab\Delta), \quad r(ab)=r(b\Delta),
\quad r(ba)=r(a\Delta) \:.
\]
In particular, we can rewrite the Traffic Equations in terms
of $r(a), r(b), r(ab),$ and $r(ba)$ only.
Observe also that: $R(a)+R(b)=R(ab)+R(ba)=1/2$.

Setting $q=1/2-p$ in \cite[(28)]{MaManta}, we
get: $\gamma_{\Delta}= -(1/2-p)R(a) -pR(b)$. In particular,
$\gamma_{\Delta}<0$. Now using \cite[(29)]{MaManta}, we get:
$\gamma= 2(1/2-p)R(b)+pR(a)$. Since $R(a)+R(b)=1/2$, we deduce that:
\[
\gamma= p + (1-4p)R(a) \:.
\]
Using the Equations in \cite[(28)]{MaManta}, we now obtain:
\[
\gamma_{\Sigma}=-2\gamma_{\Delta} = \gamma, \quad \gamma_{S_+} = 3p/2
+ (3/2-6p)R(a)\:.
\]

The only remaining point is to determine
$R(a)$. To that purpose, we need to transform the Traffic
Equations \cite[(30)]{MaManta} by switching to new unknowns.

\medskip

Recall that $R(a)=r(a)+r(a\Delta)=r(a)+ r(ba)$.
Set $R(1)= 2pr(a) + (1-2p)r(b)$ and $R(2)= r(a)+r(b)$. We have:
\[
\left[ \begin{array}{cccc}
1 & 1 & 0 & 0 \\
2p & 0 & 1-2p & 0 \\
1 & 0 & 1 & 0 \\
1 & 1 & 1 & 1
\end{array}\right] \left[ \begin{array}{c}
r(a) \\
r(ba) \\
r(b) \\
r(ab)
\end{array} \right]  = \left[ \begin{array}{c}
R(a) \\
R(1) \\
R(2) \\
1/2
\end{array} \right]\:.
\]
Assume that $p\neq 1/4$. Then, the above matrix is invertible, which
enables to write $r(a),r(b),r(ab),$ and $r(ba)$ in function of
$R(a),R(1),$ and $R(2)$. Replacing in the Traffic Equations, we get a
new set of Equations in the unknowns $R(a), R(1)$, and $R(2)$.
Select in this new set, the Equations originating from the
Equations for $r(a),
r(b)$, and $r(ab)$ in \cite[(30)]{MaManta}.
Using Maple, solve these three Equations in the unknowns $R(a), R(1)$, and
$R(2)$. We obtain that $R(a)$ is a root of the polynomial:
\[
P=2(4p-1)X^3 +
  (24p^2-18p+1)X^2+ p(-12p+7)X + p(2p-1)\:.
\]
By using that $R(2)$ should be less than 1, we obtain that for
$p<1/4$, $R(a)$ should be equal to the smallest root of $P$. This
completes the proof.
\end{proof}

\section{Proof of Proposition 4.5}\label{se-prop45}

We now prove the following statement which is Prop. 4.5 in \cite{MaManta} :

\medskip

{\bf Proposition 4.5}
{\it Assume that $\nu(a)=\nu(b)=p$,
  $\nu(a^{-1})=\nu(b^{-1})=1/2-p$, with $p\in (0,1/2)$. We have:
\begin{eqnarray*}
\gamma_{\Sigma}(p) & = & \frac{-1+ \sqrt{16p^2-8p+5}}{4}, \quad
\gamma_{S^+}(p) \ = \ \frac{4p^2+p+1-3p\sqrt{16p^2-8p+5}}{2(1-4p)} \\
\gamma_{\Delta}(p) & = & \frac{12p^2
  -5p+1-p\sqrt{16p^2-8p+5}}{2(1-4p)}\:.
\end{eqnarray*}
And eventually:
\begin{equation*}
\gamma(p)  =
\max \bigl[ 1-4p ,\frac{(1-2p)(-1-4p+\sqrt{5-8p+16p^2})}{2(1-4p)}, \frac{p
(-3+4p+\sqrt{5-8p+16p^2})}{-1+4p},-1+4p \bigr]\: .
\end{equation*} }

We give two proofs of this proposition. The first one is in
the same spirit as the proof of \cite[Prop. 4.4]{MaManta}, namely solving
the Traffic Equations \cite[(30)]{MaManta} thanks to the simplifications
provided by the symmetry $\nu(a)=\nu(b)$, $\nu(a^{-1})=\nu(b^{-1})$.
The second proof relies on the fact that, under the above symmetry,
the random process $(\widehat{X}_n)_n$ is Markovian and is a NNRW on the free
product $\Z/3\Z*\Z/3\Z$ for which the harmonic measure and the drift are known,
see \cite{MaMa}.

\medskip

{\it First proof} - We assume now that
$p=q=\nu(a)=\nu(b)=1/2-\nu(a^{-1})=1/2-\nu(b^{-1})$.
Let $r$ be the unique solution
in $\{x \in (\R_+^*)^{\Sigma} \mid  \sum_{u\in \Sigma} x(u) =1\}$
to the Traffic Equations \cite[(30)]{MaManta}.
Again, by symmetry, we should have:
\[
r(a)= r(b), \quad r(ab)=r(ba), \quad r(a\Delta)=r(b\Delta),
\quad r(ab\Delta)=r(ba\Delta) \:.
\]
Moreover, $S_a=S_b=1/2$, $R(a)=R(b)$ and $R(ab)=R(ba)=1/2-R(a)$.
Therefore, according to \cite[(28)]{MaManta} and \cite[(29)]{MaManta}, the various drifts are given by
\begin{equation}\label{gam-simple1}
\gamma_{\Sigma}=1-2p+(4p-1)R(a)\, , \: \gamma_{S_+}=1/2-2p+6pR(a)\, , \: \gamma_{\Delta}=2p-1/2-2pR(a)\, ,
\end{equation}
\begin{equation}\label{gam-simple2}
\text{for}\:p\geq 1/4\, ,\:\gamma = \begin{cases} 4p -1 & \text{if } \gamma_{\Delta} \geq 0
  \\
                       4pR(a) & \text{if } \gamma_{\Delta}
  < 0 \end{cases}\:\: ,\:\text{and for}\:p\leq1/4\, , \:\gamma(p)=\gamma(1/2-p)\:.
\end{equation}
The computation of $R(a)$ from the Traffic Equations \cite[(30)]{MaManta} in this case
turns out to be simpler than in the proof of \cite[Prop. 4.4]{MaManta}. Add the first
and the fifth equations in \cite[(30)]{MaManta}. Then $R(a)$ is a root of the
following polynomial of degree two:
\[
P\:=\:4(1-4p)X^2\:+\:2(4p-3)X\:+\:1\: .
\]
For $p\in (0,1/2)$, the only root of $P$ which lies between $0$ and $1$
is
\[
R(a)\:=\:\frac{3-4p-\sqrt{16p^2-8p+5}}{4(1-4p)}\:.
\]
Substituting in \eref{gam-simple1} and \eref{gam-simple2} completes the proof.

\medskip

\begin{figure}[ht]
\[ \epsfxsize=440pt \epsfbox{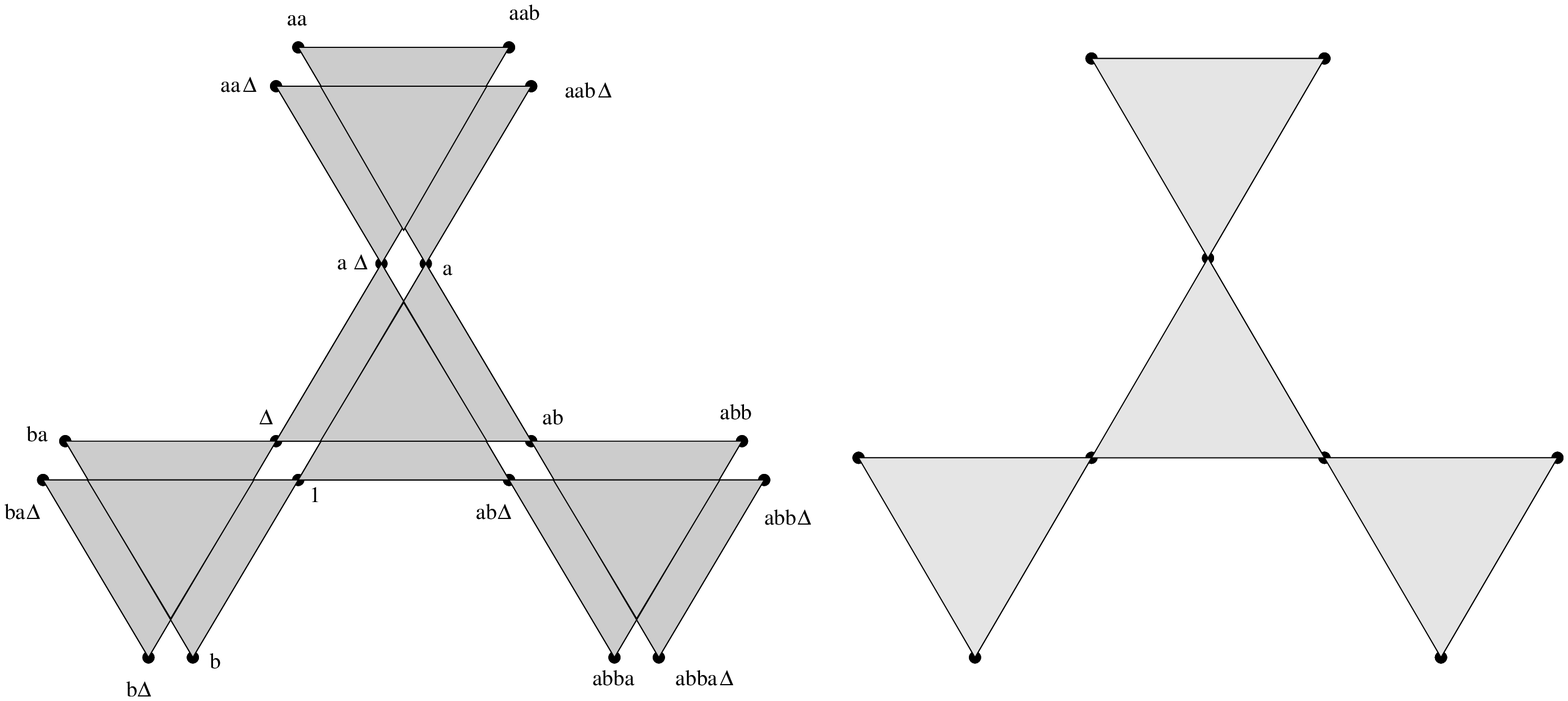} \]
\caption{The graphs $\cX(B_3/Z,S)$ (left) and $\cS(C_3,S)$ (right).}\label{fi-b3Z}
\end{figure}

{\it Second proof} - We give now another way for computing the
drift $\gamma_{\Sigma}$. Consider the subgroup $H=\{1,\Delta\}$ of
$B_3/Z$. Observe that $H$ is not a normal subgroup.
For instance, the left-class $aH=\{a,a\Delta\}$ is
different from the right-class $Ha=\{a,b\Delta\}$.
Let $C_3$ be the left-quotient of $B_3/Z$ with respect to $H$.
The elements of $C_3$ are the left-classes $\{g,g\Delta\}$.
The group $B_3$ acts on $C_3$ by left multiplication.
Denote by $\cS(C_3,S)$ the Schreier graph with respect to this action.
The set of nodes is $C_3$ and the set of
arcs is defined by:
\begin{equation}\label{schreier}
C \longrightarrow \ D \ \text{ in } \ \cS(C_3,S) \qquad \text{if}
\qquad \exists c\in C, d\in D, \
c \longrightarrow \ d \ \text{ in } \ \cX(B_3/Z,S) \:.
\end{equation}
Consider the free product $\Z/3\Z\star Z/3\Z$ with the canonical
set of generators $\widetilde{S}=\{c,c^{-1},d,d^{-1}\}$.
The Cayley graph $\cX(\Z/3\Z\star Z/3\Z,\widetilde{S})$ and the
Schreier graph $\cS(C_3,S)$ are isomorphic as unlabelled graphs,
see Figure \ref{fi-b3Z}. Below, we identify the element $\{g,g\Delta\}$
with $\widehat{g}$ where $\widehat{g}\Delta^{0/1}$ is the Garside normal form
of $g\in B_3/Z$.

Let $\nu$ be a probability distribution on $S=\{a,a^{-1},b,b^{-1}\}$.
Let $(X_n)_n$ be a realization of $(B_3,\nu)$. View $(X_n)_n$ as a
random walk on $\cX(B_3,S)$. The sequence $(p(X_n))_n$ is a realization
of the random walk $(B_3/Z,\mu)$ where $\mu=\nu\circ p^{-1}$.
It is a random walk on $\cX(B_3/Z,S)$. Recall that $\widehat{X}_n\Delta^{k_n}$
denotes the Garside normal form of $X_n$. With the identification above,
the random process $(\widehat{X}_n)_n$ evolves on $\cS(C_3,S)$ and is the
process induced by $(p(X_n))_n$. This random process $(\widehat{X}_n)_n$
is {\it a priori} not Markovian, but it is Markovian if
\[
\forall C,D \in C_3, \forall c_1,c_2\in C,\quad
P\{ p(X_{n+1}) \in D \mid p(X_n) = c_1\} = P\{ p(X_{n+1}) \in D \mid
p(X_n) = c_2\} \:.
\]
Clearly, this holds if and only if $\nu(a)=\nu(b)=p$ and $ \nu(a^{-1})=\nu(b^{-1})=1/2-p$,
which is precisely what we assume.

\medskip

Define $\cS(C_3,\mu)$ as the graph $\cS(C_3,S)$ with labels in
$[0,1]$ such that:
\[
u \stackrel{\mu(a)}{\longrightarrow} v \ \text{ in } \ \cS(C_3,\mu)
\qquad \text{if} \qquad u \stackrel{a}{\longrightarrow} v \ \text{
  in } \ \cS(C_3,S) \:.
\]
Consider the group $\Z/3\Z\star \Z/3\Z$ with generators
$\widetilde{S}=\{c,c^{-1},d,d^{-1}\}$, and the probability measure
$\tilde{\mu}$ defined on $\widetilde{S}$ by $\tilde{\mu}(c)=\tilde{\mu}(d)=\nu(a)=\nu(b),$ and
$\tilde{\mu}(c^{-1})=\tilde{\mu}(d^{-1})=\nu(a^{-1})=\nu(b^{-1})$.
Define the labelled graph $\cX(\Z/3\Z\star \Z/3\Z,\tilde{\mu})$ accordingly.
Then $\cS(C_3,\mu)$ and $\cX(\Z/3\Z\star \Z/3\Z,\tilde{\mu})$ are isomorphic
as labelled graphs. In particular, $(\widehat{X}_n)_n$ behaves like the
random walk $(\Z/3\Z\star \Z/3\Z,\tilde{\mu})$. Let us be more precise.

\medskip

For $\xi\in\widetilde{S}$, we set $\imath(\xi)=1$ if $\xi = c$ or $c^{-1}$
and $\imath(\xi)=2$ if $\xi = d$ or $d^{-1}$. Recall that $T=\{a,b,ab,ba\}\subset B_3$.
Define
\begin{eqnarray*}
L = \{\xi_0\cdots\xi_l \in T^*\mid
\text{Last}(\xi_{i})\:=\:\text{First}(\xi_{i+1})\}, && L^{\infty} =  \{\xi_0\cdots \in T^{\N}\mid
\text{Last}(\xi_{i})\:=\:\text{First}(\xi_{i+1})\} \\
\widetilde{L} =  \{\xi_0\xi_1\cdots \xi_l \in
\widetilde{S}^* \mid
\imath(\xi_{i})\: \neq \:\imath(\xi_{i+1})\}, && \widetilde{L}^{\infty}  = \{\xi_0\xi_1\cdots \in
\widetilde{S}^{\N} \mid
\imath(\xi_{i})\: \neq \:\imath(\xi_{i+1})\} \:.
\end{eqnarray*}

Set $\lim_n \widehat{X}_n = \widehat{X}_{\infty} = \widehat{x}_0 \widehat{x}_1 \widehat{x}_2\cdots$,
with $\widehat{x}_i\in T$. Let
$\mu^{\infty}$ be the law of $\widehat{X}_{\infty}$; it is a measure
on $L^{\infty}$. Let $\tilde{\mu}^{\infty}$
be the harmonic measure of $(\Z/3\Z\star \Z/3\Z,\tilde{\mu})$; it is a measure
on $\widetilde{L}^{\infty}$.
Let $\widetilde{r}$ be the unique
solution in $\mathring{\cB}$ of the Traffic Equations of
$(\Z/3\Z\star \Z/3\Z,\tilde{\mu})$, see \cite[\S 4.3]{MaMa}.
Define
\[
r(a)\:=\:r(b)\:=\:\tilde{r}(c)\:=\: \tilde{r}(d)\: ,\:
r(ab)\:=\:r(ba)\:=\:\tilde{r}(c^2)\:=\: \tilde{r}(d^2)\: ,
\]
\[
r(1)\:=\:r(a)+r(b)\: ,\:\text{and}\quad\:r(2)\:=\:r(ab)+r(ba) \:.
\]
Then, we have, $\forall u_1\cdots u_l \in L$,
\[
\mu^{\infty} (u_1u_2u_3\cdots u_l T^{\N}) =  \tilde{\mu}^{\infty}
(c^{|u_1|_{S}}d^{|u_2|_{S}}c^{|u_3|_{S}} \cdots
\widetilde{\Sigma}^{\N} )
=  (1/2) r(|u_1|_{S})\cdots r(|u_l|_{S})\:.
\]
In particular
\begin{equation}\label{eq-ouf}
\gamma_{\Sigma} = \lim_{n\rightarrow \infty} \frac{|\widehat{X}_n|_T}{n} =
\lim_{n\rightarrow \infty} \frac{|W_n|_{\widetilde{S}}}{n}=
pr(1) + \bigl(\frac{1}{2}-p\bigr)r(2) \:,
\end{equation}
where $(W_n)_n$ is a realization of the random walk $(\Z/3\Z\star \Z/3\Z,\tilde{\mu})$.
See \cite[Corollary 3.6]{MaMa} for the third equality in Eq. \eref{eq-ouf}.

\medskip

The vector $\tilde{r}$ is computed in \cite[\S 4.3]{MaMa}:
\begin{equation}\label{eq-z3z3a=b}
\tilde{r}(c)=\tilde{r}(d)=\frac{4p-3+\sqrt{16p^2-8p+5}}{4(4p-1)}, \quad
\tilde{r}(c^2)=\tilde{r}(d^2)=\frac{4p+1-\sqrt{16p^2-8p+5}}{4(4p-1)}\:.
\end{equation}
We get $\gamma_{\Sigma}\:=\:(-1 + \sqrt{16p^2-8p+5})/4$. We deduce
the other drifts easily.

\section{Proof of Proposition 5.4}\label{se-prop54}

We now prove the following statement which is Prop. 5.4 in \cite{MaManta} :

\medskip

{\bf Proposition 5.4}
{\it Consider the simple random walk $(A_k,\nu)$ with
$\nu(a)=\nu(a^{-1})=\nu(b)=\nu(b^{-1})=1/4$. The drifts $\gamma_{\Sigma}$ and
$\gamma_{\Delta}$ are given by
\begin{equation}
\gamma_{\Sigma} = \frac{1-x_k}{2}, \qquad \gamma_{\Delta}=
-\frac{1-x_k}{4}\:.
\end{equation}
Let $\gamma$ be the drift of the length with respect to the natural
generators $\{a,b,a^{-1},b^{-1}\}$. We have
\begin{equation}
\gamma = \begin{cases}
(1-x_k) \bigl[\ \sum_{i=1}^{j-1} i F_i(x_k) + (j/2) F_{j}(x_k)\ \bigr]
  & \text{if } \ k=2j \\
(1-x_k) \bigl[ \ \sum_{i=1}^{j} i F_i(x_k) \ \bigr]& \text{if } \ k=2j+1
\end{cases}\:.
\end{equation} }

\begin{proof}
Consider the group $\Z/k\Z\star \Z/k\Z$
with generators $\widetilde{S}=\{c,c^{-1},d,d^{-1}\}$. Let $(X_n)_n$ be a
realization of the simple random walk on $A_k$, and recall that
we write $\widehat{X}_n\Delta^{k_n}$ for the Garside normal form of $X_n$.
Denote by $\mu$ and $\tilde{\mu}$ the uniform probability measures on
$S$ and $\widetilde{S}$ respectively.

\par

If $k$ is even, then the unlabelled Cayley graph
$\cX(A_k/Z,S)$ is isomorphic to the unlabelled Cayley graph
$\cX(\Z/k\Z\star Z/k\Z,\widetilde{S})$ (see Figure \ref{fi-A4}) so the
simple random walks $(A_k/Z,\mu)$ and $(\Z/k\Z\star \Z/k\Z,\tilde{\mu})$
are isomorphic.

\begin{figure}[ht]
\[ \epsfxsize=400pt \epsfbox{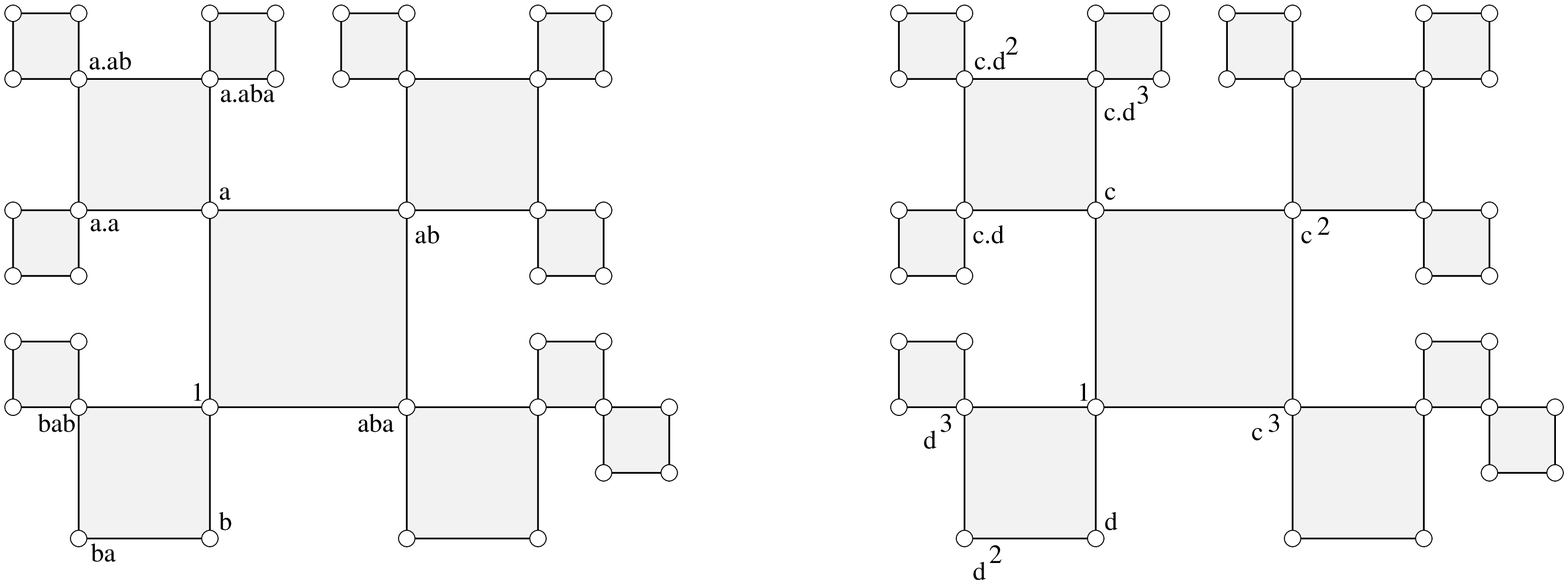} \]
\caption{The Cayley graphs $\cX(A_4/Z,S)$ (left) and $\cX(\Z/4\Z\star
  \Z/4\Z,\widetilde{S})$ (right).}\label{fi-A4}
\end{figure}

For $k$ odd, the structure of $\cX(A_k/Z,S)$ is more
twisted. This is illustrated in Figure \ref{fi-A5b}
(see also Figure \ref{fi-b3Z} (left)).

\begin{figure}[ht]
\[ \epsfxsize=180pt \epsfbox{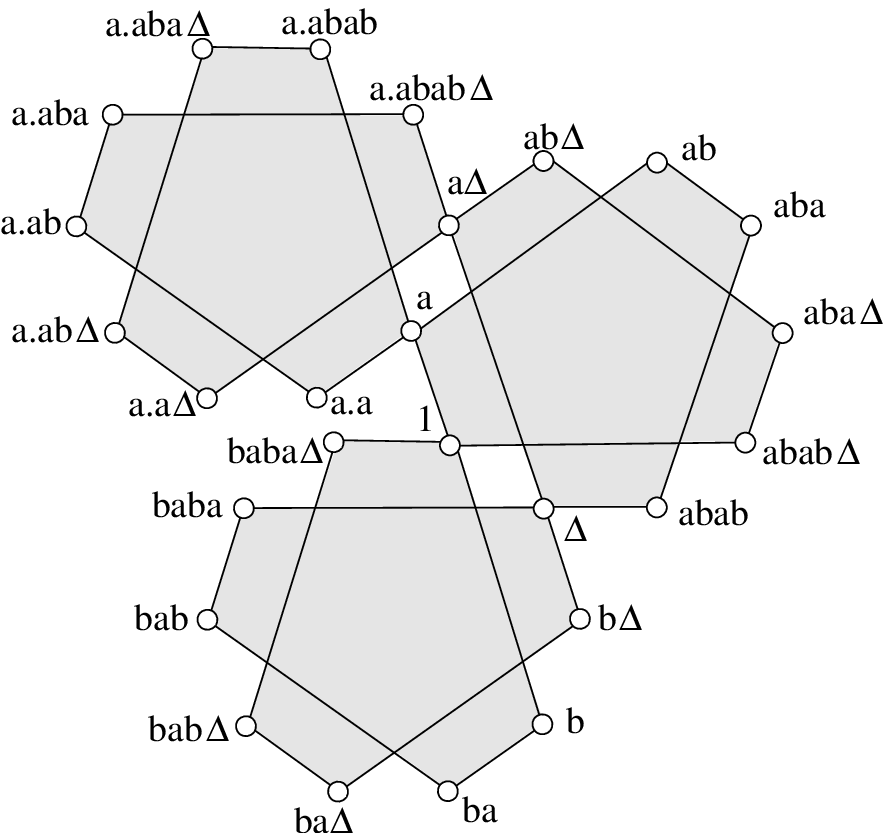} \]
\caption{The Cayley graph $\cX(A_5/Z,S)$.}\label{fi-A5b}
\end{figure}

Let $C_k$ be the left-quotient of
$A_k/Z$ by the (non-normal) subgroup $H=\{1,\Delta\}$, and define
the Schreier graph $\cS(C_k,S)$ as in \eref{schreier}.
Clearly $\cS(C_k,S)$ is isomorphic to the Cayley graph
$\cX(\Z/k\Z\star Z/k\Z,\widetilde{S})$
(as unlabelled graphs). This is illustrated in Figure \ref{fi-A5ter}.

\begin{figure}[ht]
\[ \epsfxsize=440pt \epsfbox{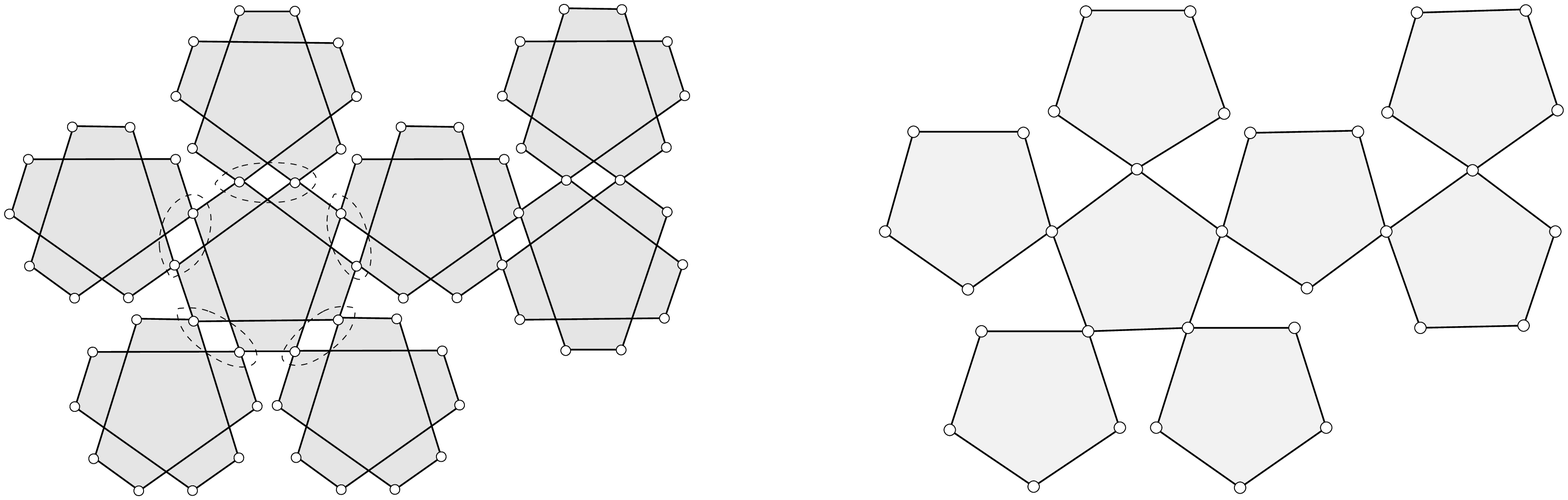} \]
\caption{The graphs $\cX(A_5/Z,S)$ (left) and $\cS(C_5,S)$ (right).}\label{fi-A5ter}
\end{figure}

Therefore the simple random walks on the two graphs $\cS(C_k,S)$ and
$\cX(\Z/k\Z\star Z/k\Z,\widetilde{S})$ are isomorphic.
\par
In both cases, the Markovian random process $\widehat{X}_n$ behaves
like the simple random walk $(\Z/k\Z\star \Z/k\Z,\tilde{\mu})$.
Adapting the end of the second proof of \cite[Prop. 4.5]{MaManta}
and using the results of \cite[\S 4.4]{MaMa} leads to
$\gamma_{\Sigma}$, hence to the other drifts.
\end{proof}

\section{Extensions}\label{se-exte}

One can retrieve from the harmonic measure $\mu^{\infty}$, other
quantities of interest for the random walk $(B_3/Z,\mu)$: (a) the entropy,
(b) the minimal positive harmonic functions, or (c) the Green function.

Consider for instance the entropy $h$. For NNRW on free products of finite
groups, or on 0-automatic pairs, a formula was available for $h$ as a
simple function of $r$, the unique solution to the Traffic
Equations, see \cite{mair04,MaMa}. Here the situation is more complex
and $h$ can only be expressed as a limit of functions of $r$. However
$h$ can be computed with an arbitrary prescribed precision.

In the three above cases (a), (b), and (c), the key is to determine the Radon-Nikodym derivatives
$du\circledast \mu^{\infty}/d\mu^{\infty} (\cdot )$ for $u\in \cG$.

\medskip

Fix $u=u_1\cdots u_k \in \cG$ and $\xi=\xi_1\xi_2\cdots \in
\cG^{\infty}$. We have, using \cite[(21)]{MaManta},
\begin{equation}\label{eq-limit}
\frac{du\circledast \mu^{\infty}}{d\mu^{\infty}} (\xi) = \lim_n
\frac{\mu^{\infty}(u^{-1}\circledast \xi_1\cdots \xi_n
  T^{\N})}{\mu^{\infty}(\xi_1\cdots \xi_n T^{\N})} = \lim_n
\frac{\alpha\cM(w_1\cdots w_{\ell-1})\beta(w_{\ell})}{\alpha
  \cM(\xi_1\cdots \xi_{n-1})\beta(\xi_n)} \:,
\end{equation}
where $w_1\cdots w_{\ell}= u^{-1}\circledast \xi_1\cdots \xi_n$.
Below, we justify the existence of the limit in \eref{eq-limit}, and in doing so we show
how to control the error made when replacing the limit by the value
computed for a given $n$.

\medskip

It follows from the definition in \cite[(9)]{MaManta} that we have either:
\[
u^{-1}\circledast \xi_1\cdots \xi_n = v\cdot \xi_l \cdot \ldots \cdot \xi_n,
\quad \text{or} \quad u^{-1}\circledast \xi_1\cdots \xi_n = v\cdot \iota(\xi_l)
\cdot \ldots \cdot \iota(\xi_n)\:,
\]
for some $v \in \cG$ and $l\in \N^*$ which do not depend on $n$, for
$n$ large
enough.
In the first case, respectively the second one, we have:
\begin{eqnarray}
\frac{du\circledast \mu^{\infty}}{d\mu^{\infty}} (\xi) &  = & \lim_n \
  \frac{\alpha\cM(v\xi_{\ell}\cdots \xi_{n-1})\beta(\xi_n)}{\alpha
  \cM(\xi_1\cdots \xi_{n-1})\beta(\xi_n)} \label{eq-radon1} \\
\text{resp. } \ \frac{du\circledast \mu^{\infty}}{d\mu^{\infty}} (\xi)
  & = & \lim_n \
  \frac{\alpha\cM(v\iota(\xi_{\ell})\cdots \iota(\xi_{n-1}))\beta(\iota(\xi_n))}{\alpha
  \cM(\xi_1\cdots \xi_{n-1})\beta(\xi_n)}\:. \label{eq-radon2}
\end{eqnarray}
Now, the $\R_+$-automaton $(\alpha,\cM,\beta)$ of \cite[Figure 9]{MaManta} has several remarkable
properties.
Define $\iota(\alpha)=[0,1,0,1]\in
\R_+^{1\times Q}$. Observe that $\iota(\alpha)_i=
\alpha_{5-i}$. Observe also that:
$\forall u \in \Sigma, \forall i,j, \ \cM(u)_{ij}=
\cM(\iota(u))_{5-i,5-j}$.
It implies that we have:
\begin{equation}\label{eq-sym}
\forall u=u_1\cdots u_k \in T^*, \qquad \alpha\cM(u_1\cdots u_{k-1})\beta(u_k)=
\iota(\alpha)\cM(\iota(u_1)\cdots \iota(u_{k-1}))\beta(\iota(u_k))\:.
\end{equation}
Extend by morphism the map $\iota$ defined in \cite[(10)]{MaManta} to
$\iota: \Sigma^*\longrightarrow \Sigma^*$.
The identity in \eref{eq-sym} enables to rewrite \eref{eq-radon2} as:
\[
\frac{du\circledast \mu^{\infty}}{d\mu^{\infty}} (\xi)
 = \lim_n \
  \frac{\iota(\alpha)\cM(\iota(v)\xi_{\ell}\cdots \xi_{n-1})\beta(\xi_n)}{\alpha
  \cM(\xi_1\cdots \xi_{n-1})\beta(\xi_n)}\:.
\]
To summarize, in all cases, for $n$ large enough, there exist $l\in
\N^*$ and $\alpha_1 \in \R_+^{1\times Q}$ such that:
\begin{equation}\label{eq-birkhoff}
\frac{du\circledast \mu^{\infty}}{d\mu^{\infty}} (\xi)
 = \lim_n \ \frac{\alpha_1 \cM(\xi_{\ell}\cdots
 \xi_{n-1})\beta(\xi_n)}{\alpha_2 \cM(\xi_{\ell}\cdots
 \xi_{n-1})\beta(\xi_n)} \:,
\end{equation}
with $\alpha_2=\alpha \cM(\xi_1\cdots \xi_{\ell-1})$.
The limit in \eref{eq-birkhoff} exists as a consequence of general
results on inhomogeneous products of
non-negative matrices, see for instance \cite[Chapter 3]{sene}. To be
more precise, and to evaluate the speed of convergence, it is
convenient to slightly rewrite \eref{eq-birkhoff}, in order to have
matrices with positive entries.

Consider the automaton $(\widetilde{\alpha}, \widetilde{\cM},
\widetilde{\beta})$ defined as follows.
Set $\widetilde{\alpha} = [1,1]$; for $u\in T$,
set $\widetilde{\beta}(u) =
[R(u),R(u)]^T$;
and let the morphism $\widetilde{\cM}:T^* \rightarrow \R_+^{2\times 2}$
be defined by:
\[
\widetilde{\cM}(a) =  \left[ \begin{array}{cc} q(a) & q(a\Delta) \\ q(b\Delta) &
    q(b) \end{array}\right] \:, \quad \widetilde{\cM}(b)  = \ \left[
    \begin{array}{cc} q(b) & q(b\Delta) \\ q(a\Delta) &
    q(a) \end{array}\right]
\]
\[
\widetilde{\cM}(ab)  = \ \left[ \begin{array}{cc} q(ab) & q(ab\Delta) \\ q(ba\Delta) &
    q(ba) \end{array}\right] \:, \quad \widetilde{\cM}(ba)  = \ \left[
    \begin{array}{cc} q(ba) & q(ba\Delta) \\ q(ab\Delta) &
    q(ab) \end{array}\right] \:.
\]
It is easily checked that, {\em on} $\cG$, the two automata $(\widetilde{\alpha}, \widetilde{\cM},
\widetilde{\beta})$ and $(\alpha,\cM,\beta)$ coincide. That is:
\begin{equation}\label{eq-coincide}
\forall u=u_1\cdots u_k  \in \cG, \quad \widetilde{\alpha}\widetilde{\cM}(u_1\cdots
u_{k-1})\widetilde{\beta}(u_k) = \alpha\cM(u_1\cdots
u_{k-1})\beta(u_k) \:.
\end{equation}
Observe that the above identity does not hold on $T^*\moins \cG$. (In
fact, the automaton  $(\alpha,\cM,\beta)$ is the tensor product of $(\widetilde{\alpha}, \widetilde{\cM},
\widetilde{\beta})$ with a 2-state automaton recognizing $\cG\cap T^*$.)

\medskip

Define:
\begin{equation*}
\delta = \min_{u\in \Sigma} \frac{ \min_{ij}
  \widetilde{\cM}(u)_{ij}}{\max_{ij} \widetilde{\cM}(u)_{ij}}, \quad K =
  \max_{u\in \Sigma} \max_{ij} \Bigl[ \max_k
  \frac{\widetilde{\cM}(u)_{ik}}{\widetilde{\cM}(u)_{jk}} - \min_k
  \frac{\widetilde{\cM}(u)_{ik}}{\widetilde{\cM}(u)_{jk}} \Bigr] \:.
\end{equation*}

Observe that $0<\delta <1$. For $x=x_1x_2\cdots \in T^{\N}$, set
\[
c(x)=\bigl[ 1,\lim_n \widetilde{\cM}(x_1\cdots
    x_n)_{2k}/\widetilde{\cM}(x_1\cdots x_n)_{1k} \bigl]^T\:,
\]
where the limit does not depend on $k\in \{1,2\}$.
Using \cite[Exercice 3.9]{sene}, we get, for $k\in \{1,2\}$,
\begin{equation}\label{eq-seneta}
\Bigl| \ \frac{\widetilde{\cM}(x_1\cdots
    x_n)_{2k}}{\widetilde{\cM}(x_1\cdots x_n)_{1k}} -c(x)_2 \ \Bigr| \leq
    K(1-\delta^2)^{n-1} \:.
\end{equation}
Now let us go back to \eref{eq-birkhoff}. Set $\eta=\eta_1\eta_2\cdots
= \xi_{\ell}\xi_{\ell+1}\cdots$. For $x=(x_1,x_2)\in \R^2,$
set $\norm{x} =|x_1|+|x_2|$; for $x=(x_1,x_2), y=(y_1,y_2)\in \R^2,$
set $\langle x,y\rangle = x_1y_1+ x_2y_2$. We have:
\[
\frac{du\circledast \mu^{\infty}}{d\mu^{\infty}} (\xi) = \lim_n
\frac{\widetilde{\alpha}_1\widetilde{\cM}(\xi_{\ell}\cdots \xi_n)
  \widetilde{\beta}(\xi_{n+1})}{\widetilde{\alpha}_2\widetilde{\cM}(\xi_{\ell}\cdots \xi_n)
  \widetilde{\beta}(\xi_{n+1})} = \frac{ \langle \widetilde{\alpha}_1, c(\eta)
  \rangle}{\langle \widetilde{\alpha}_2 , c(\eta)\rangle}\:.
\]
Let us fix $n$ larger than
$\ell$. Set $\varepsilon = K(1-\delta^2)^{n-\ell-1}$.
Using \eref{eq-seneta}, we easily get:
\begin{equation*}
\frac{ \langle \widetilde{\alpha}_1, c(\eta)
  \rangle - \varepsilon \norm{\widetilde{\alpha}_1}}{ \langle \widetilde{\alpha}_2, c(\eta)
  \rangle + \varepsilon \norm{\widetilde{\alpha}_2}} \ \leq \
  \frac{\widetilde{\alpha}_1\widetilde{\cM}(\xi_{\ell}\cdots \xi_n)
  \widetilde{\beta}(\xi_{n+1})}{\widetilde{\alpha}_2\widetilde{\cM}(\xi_{\ell}\cdots \xi_n)
  \widetilde{\beta}(\xi_{n+1})} \ \leq \ \frac{ \langle \widetilde{\alpha}_1, c(\eta)
  \rangle + \varepsilon \norm{\widetilde{\alpha}_1}}{ \langle \widetilde{\alpha}_2, c(\eta)
  \rangle - \varepsilon \norm{\widetilde{\alpha}_2}} \:.
\end{equation*}
The above inequalities provide a sharp control on the error made when
replacing the limit in \eref{eq-birkhoff} by the value computed for a fixed $n$.

\paragraph{Entropy} $ $ \medskip

The {\em entropy} of a probability measure $\mu$ with finite support $S$ is
defined by $H(\mu) = -\sum_{x\in S} \mu(x) \log [\mu(x) ]$.
Consider a random walk $(G,\mu)$, defined as in \cite[Section 2]{MaManta}.
Let $(X_n)_n$ be a realization of the random walk.
The {\em entropy of $(G,\mu)$}, introduced by Avez \cite{avez}, is
\begin{equation}\label{eq-avez}
h = \lim_n \frac{H(\mu^{*n})}{n}= \lim_n -\frac{1}{n}
\log\mu^{*n}(X_n) \:,
\end{equation}
$a.s.$ and in $L^p$, for all $1\leq p <\infty$. The existence of the
limits as well as their equality follow from Kingman's subadditive
ergodic theorem~\cite{avez,derr}.

\medskip

Consider the random walk $(B_3/Z,\mu)$.
Recall that
$(\cG^{\infty},\mu^{\infty})$ is the Poisson boundary of
$(B_3/Z,\mu)$. Then, we have, see \cite[Theorem 3.1]{KaVe}:
\begin{equation}\label{eq-kave}
h =  - \sum_{u\in \Sigma} \mu(u) \int \log \bigl[
  \frac{du^{-1}\circledast\mu^{\infty}}{d\mu^{\infty}}(\xi)\bigr]
  d\mu^{\infty}(\xi) \:.
\end{equation}
Using \eref{eq-kave} together with the above remarks on how to
approximate the Radon-Nikodym derivatives, one can derive an algorithm
to compute $h$ with an arbitrarily prescribed precision.

\paragraph{Minimal positive harmonic functions} $ $ \medskip

Consider a random walk $(G,\mu)$ defined as in \cite[Section 2]{MaManta}.
A {\em positive harmonic} function is a function
$f: G \rightarrow \R_+$ such that: $\forall u \in G, \
 \sum_{a\in \Sigma} f(u\ast a)\mu(a) = f(u)$. A positive harmonic function $f$ is {\em
   minimal} if $f(1)=1$ and if for any positive harmonic function $g$ such that $f\geq
 g$, there exists $c\in \R_+$ such that $f=cg$.

\medskip

Consider the random walk $(B_3/Z,\mu)$ as above. Recall that $\cG^{\infty}$ is
the minimal Martin boundary of $(B_3/Z,\mu)$. So the set of minimal
positive harmonic function is precisely given by $\{K_{\xi}, \xi \in
\cG^{\infty}\}$ with
\begin{equation}\label{eq-harmonic}
K_{\xi}: B_3/Z \rightarrow \R_+, \quad K_{\xi}(g) = \frac{d
  \phi(g)\circledast \mu^{\infty}}{d\mu^{\infty}}(\xi) \:.
\end{equation}

\paragraph{Green function} $ $ \medskip

Consider a transient random walk $(G,\mu)$ defined as in \cite[Section 2]{MaManta} and a realization $(X_n)_n$.
Define, $Q:G \rightarrow [0,1],$
\begin{equation}\label{eq-ever}
Q(g)=  P \{ \exists n \geq 1 \mid X_n = g\}\:,
\end{equation}
the probability of ever reaching $g$. The {\em Green
function} is the map $\Gamma: G \rightarrow \R_+, \
\Gamma(g)=\sum_{i=0}^{\infty} \mu^{*i}(g)$. Observe that:
\begin{equation}\label{eq-green}
\Gamma(\cdot)= \Gamma(1)Q(\cdot)\:.
\end{equation}

Consider now the random walk $(B_3/Z,\mu)$.
Computing $Q$ is more involved than for
free groups \cite{DyMa} or for zero-automatic pairs
\cite{mair04}. However, the spirit remains the same: there is a close
link between $Q$ and $r$, the unique solution to the Traffic
Equations.

For convenience, we view $Q$ and $\Gamma$ as $Q:\cG \rightarrow \R_+$
and $\Gamma:\cG \rightarrow \R_+$. Define, for $u\in \cG$, the auxiliary quantity:
\begin{eqnarray}\label{eq-aux}
\widehat{q}(u) & = & P \{ \exists n > 1 \mid Y_n = u \ \text{ and
} \ \forall 1 < m < n, Y_m \neq u\Delta \}\:.
\end{eqnarray}

We have the following:

\begin{proposition}\label{pr-ever}
Let $r \in
\{ x \in (\R_+^*)^{\Sigma} \mid  \sum_{u\in \Sigma} x(u) =1\}$ be the
unique solution to the Traffic Equations \cite[(19)]{MaManta} of
$(B_3/Z,\mu)$. For $u\in \Sigma$, set $q(u)= r(u)/r(\next(u))$.
We have:
\begin{equation}\label{eq-q1}
\forall u\in \Sigma, \quad \widehat{q}(u) = \frac{r(u)}{r(\next(u))} =q(u)\:.
\end{equation}
Besides:
\begin{equation}\label{eq-q2}
\widehat{q}(1) = \sum_{u\in \Sigma} \mu(u) q(u^{-1}), \quad
\widehat{q}(\Delta) = \sum_{u\in \Sigma} \mu(u) q(u^{-1}\Delta)\:.
\end{equation}
The probabilities of ever reaching an element are given by: $Q(\Delta)
= \widehat{q}(\Delta)/(1-\widehat{q}(1))$, and, $\forall v=v_1\cdots v_k
\in \cG \moins \{1,\Delta\}$:
\begin{equation}\label{eq-q3}
Q(v) =  \sum_{u_1\cdots u_k \in \psi^{-1}(v)} q(u_1)\cdots
q(u_{k-1}) \bigl[ q(u_{k}) +
q(u_{k}\Delta)Q(\Delta)
\bigr] \:.
\end{equation}
The Green function is determined by \eref{eq-green} and:
\begin{equation}\label{eq-green2}
\Gamma(1) = \frac{1-\widehat{q}(1)}{(1-\widehat{q}(1))^2 -
  \widehat{q}(\Delta)^2}\:.
\end{equation}
\end{proposition}

(See \cite[Eq. (13)]{MaManta} for the definition of $\psi$)
\begin{proof}
It follows from the shape of the Cayley graph $\cX(B_3/Z,\Sigma)$
that:
\[
\widehat{q}(1) = \sum_{u\in \Sigma} \mu(u) \widehat{q}(u^{-1}), \quad
\widehat{q}(\Delta) = \sum_{u\in \Sigma} \mu(u) \widehat{q}(u^{-1}\Delta), \quad Q(\Delta)
= \widehat{q}(\Delta)/(1-\widehat{q}(1))\:.
\]
Similarly, \eref{eq-q3} holds with $\widehat{q}(\cdot)$ in place
of $q(\cdot)$. Let us prove \eref{eq-green2}. Define
$\bar{Q}(1)=P\{\exists n > 1 \mid Y_n =1\}$. Clearly, $\bar{Q}(1)=
\widehat{q}(1) + \widehat{q}(\Delta) Q(\Delta)$ and $\Gamma(1)= 1+
\bar{Q}(1) \Gamma(1)$. The expression in \eref{eq-green2} follows.
Therefore, the only
point that remains to be proved is: $\forall u\in
\Sigma, \widehat{q}(u)=q(u)$.

\medskip

For $u \in \Sigma$, define $F_u =\{ v\in \Sigma \mid \first(v)=\first(u)\}$.
By considering the random walk after one move, we get
that $(\widehat{q}(u))_{u\in \Sigma}$ is a solution to
the following
equations over the indeterminates $(y(u))_{u\in \Sigma}$:
\begin{equation}\label{eq-widehatq}
y(u) = \mu(u) + \sum_{v\in F_u \moins \{u,u\Delta\}} \mu(v)y(v^{-1}u) +
\sum_{v\in \Sigma\moins F_u} \mu(v) \bigl[ y(v^{-1})y(u)+
y(v^{-1}\Delta)y(\iota(u)\Delta) \bigr] \:.
\end{equation}
Starting from the Traffic Equations \cite[(19)]{MaManta}, and dividing
by $x(\next(u))$, we obtain precisely the Equations
\eref{eq-widehatq}.
Let $r \in
\{ x \in (\R_+^*)^{\Sigma} \mid  \sum_{u\in \Sigma} x(u) =1\}$ be the
unique solution to the Traffic Equations, and set $q(u) =
r(u)/r(\next(u))$ for $u\in \Sigma$.
We deduce from the above that $(q(u))_{u \in \Sigma}$ is a solution
to the Equations \eref{eq-widehatq}.

\medskip

We cannot directly conclude that $\widehat{q}(u)=q(u)$. Indeed the
Equations \eref{eq-widehatq} do not characterize
 $(\widehat{q}(u))_{u\in \Sigma}$. They
have in general several solutions in
$(0,1)^{\Sigma}$
(if $\mu$ is uniform over $\Sigma$, then the two constant functions
$y=1/2$ and $y=1/4$ satisfy \eref{eq-widehatq}).
This is in contrast with the situation
for 0-automatic pairs where the analogs of Equations \eref{eq-widehatq}
have a unique solution \cite[Lemma 4.7]{mair04}.

\medskip

Recall that the minimal
Martin boundary coincides with the Martin boundary and is
$\cG^{\infty}$. Hence, the minimal positive
harmonic functions, given in \eref{eq-harmonic}, can also be described as:
$\forall \xi=\xi_1\cdots \in \cG^{\infty}$,
\begin{equation}\label{eq-harmonic2}
K_{\xi} : B_3/Z \rightarrow \R_+, \quad K_{\xi}(g) = \lim_n
  \frac{\Gamma(\phi(g^{-1})\circledast \xi_1\cdots \xi_n)}{
  \Gamma(\xi_1\cdots \xi_n)}= \lim_n \frac{Q(\phi(g^{-1})\circledast \xi_1\cdots
  \xi_n)}{Q(\xi_1\cdots \xi_n)} \:.
\end{equation}
Juxtaposing \eref{eq-harmonic} and \eref{eq-harmonic2}, and using
 \cite[(20)]{MaManta} and \eref{eq-q3} (with $\widehat{q}(\cdot)$ replacing  $q(\cdot)$), we
get: for all $u\in B_3/Z$ and $\xi=\xi_1\xi_2\cdots \in \cG^{\infty}$,
\begin{eqnarray*}
K_{\xi}(u) & = & \lim_n \ \frac{\sum_{v_1\cdots v_{\ell} \in
    \psi^{-1}(\phi(u^{-1})\circledast \xi_1\cdots \xi_n)} q(v_1)\cdots
    q(v_{\ell-1}) R(v_{\ell})}{\sum_{v_1\cdots v_n \in
    \psi^{-1}(\xi_1\cdots \xi_n)} q(v_1)\cdots
    q(v_{n-1}) R(v_{n})}  \\
& = & \lim_n \ \frac{\sum_{v_1\cdots v_{\ell} \in
    \psi^{-1}(\phi(u^{-1})\circledast \xi_1\cdots \xi_n)} \widehat{q}(v_1)\cdots
    \widehat{q}(v_{\ell-1}) Q(v_{\ell})}{\sum_{v_1\cdots v_n \in
    \psi^{-1}(\xi_1\cdots \xi_n)} \widehat{q}(v_1)\cdots
    \widehat{q}(v_{n-1}) Q(v_{n})} \:.
\end{eqnarray*}
By choosing appropriately the values of $u$, we deduce easily that
this implies $q(\cdot )=\widehat{q}(\cdot )$.
\end{proof}

\paragraph{Central Limit Theorem} $ $ \medskip

Recall that we write $\widehat{X}_n\Delta^{k_n}$ for the Garside
normal form of $X_n$. The description of the asymptotic behavior
the quotient process $(p(X_n))$ evolving on $B_3/Z$ provides a
Central Limit Theorem for both the length $|\widehat{X}_n|_T$
and the exponent $k_n$ in the same way as in \cite{ledr00}.
\par
The statement is the following. Recall that we set
$U=p^{-1}(\Sigma)=\{a\Delta^k,ab\Delta^k,b\Delta^k,$
$ba\Delta^k,k\in\Z\}$ and $T=\{a,b,ab,ba\}\subset B_3$

\begin{proposition}\label{CLT}
Consider the random walk $(B_3,\nu)$ where $\nu$ is a probability
measure on $U$ such that $\cup_n \text{supp} (\nu^{*n})=B_3$ and
$\sum_{x\in U} e^{\lambda |x|_S} \nu(x)< \infty$
for some $\lambda > 0$. Let $(X_n)$ be a realization of the
random walk $(B_3,\nu)$ and $\widehat{X}_n\Delta^{k_n}$ the Garside
normal form of $X_n$. Then there exist
two positive numbers $\sigma_{\Sigma}$ and $\sigma_{\Delta}$ such that, for all $t\in \R$,
\begin{eqnarray}
\lim_{n\to +\infty}{P\biggl\{\frac{|\widehat{X}_n|_{T}-n\gamma_{\Sigma}}{\sigma_{\Sigma}\sqrt{n}} < t \biggr\} }
\:=\:
\lim_{n\to +\infty}{P\biggl\{\frac{k_n - n \gamma_{\Delta}}{\sigma_{\Delta}\sqrt{n}} < t \biggr\} }
\:=\:
\frac{1}{\sqrt{2\pi}}\int_{-\infty}^{t}e^{-x^2/2}\text{d}x \:\:.
\end{eqnarray}
\end{proposition}

As far as $k_n$ is concerned, starting from the proof of \cite[Lemma 4.10]{MaManta},
one can easily adapt the techniques developped in \cite[Section 4.c]{ledr00}
to prove the result. The proof of the statement for $|\widehat{X}_n|_{T}$ is analogous,
the function $\theta_{\Delta}$ has to be replaced by a
function $\theta_{\Sigma}\: :\: U^{\N} \times \cG^{\infty} \rightarrow \Z$ defined,
for all $(\omega,\xi)\in U^{\N} \times \cG^{\infty}$, by
\begin{equation*}
\theta_{\Sigma} (\omega,\xi)=
    \begin{cases} +1 & \mbox{if}\:\:\omega_{0}\ast \xi_0 \not\in \Sigma \cup \{1,\Delta\} \\
                  -1 & \mbox{if}\:\:\omega_{0}\ast \xi_0 \:=\:1\:\mbox{or}\:\Delta \\
                  0  & \mbox{if}\:\:\omega_{0}\ast \xi_0 \in \Sigma
    \end{cases}\: ,\:\:\mbox{where} \:\omega\:=\:(\omega_0,\omega_1,...)\:.
\end{equation*}
That is, $\theta_{\Sigma} (\omega,\xi)$ counts the variation of the length of $\xi$
when it is left-multiplied by $\omega_0$. Integrating $\theta_{\Sigma}$ over
$U^{\N} \times \cG^{\infty}$ leads to the drift $\gamma_{\Sigma}$, exactly  as in \cite[Lemma 4.10]{MaManta}.

\end{document}